\documentclass[12pt,twoside]{article}
\usepackage{amsmath,amssymb,amsthm}
\usepackage{epsfig,graphicx,verbatim}
\usepackage{color}
\usepackage[bookmarks, bookmarksopen=true, bookmarksnumbered=true]{hyperref}
\numberwithin{equation}{section}
\newtheorem{thm}[equation]{Theorem}

\newtheorem{defn}[equation]{Definition}
\newtheorem{conj}[equation]{Conjecture}
\newcounter{mycount}
\newenvironment{romlist}{\begin{list}{\rm(\roman{mycount})}%
   {\usecounter{mycount}\labelwidth=1cm\itemsep 0pt}}{\end{list}}
\newenvironment{numlist}{\begin{list}{\arabic{mycount}.}%
   {\usecounter{mycount}\labelwidth=1cm\itemsep 0pt}}{\end{list}}

\newcommand{\tr}{\operatorname{tr}}
\newcommand{\rad}{\operatorname{rad}}

\newcommand{\uinvnorm}{|\kern-2pt|\kern-2pt|}

\def\1{\mbox{1\hskip-.25em l}}
\newcommand{\beq}{\begin{equation}}
\newcommand{\eeq}{\end{equation}}
\newcommand{\ZZ}{\mathbb{Z}}

\theoremstyle{plain}

\theoremstyle{definition}

\theoremstyle{remark}

\def\anb{a_{n,\b}}
\def\amb{a_{m,\b}}

\def\sF{\mathcal F}

\def\sH{\mathcal H}

\def\qq{\qquad}

\def\a{\alpha}
\def\b{\beta}
\def\de{\delta}
\def\lam{\lambda}
\def\lamc{\lam_{\mathrm c}}
\def\vlamc{\vec\lam_{\mathrm c}}
\def\th{\theta}
\def\s{\sigma}

\def\rc{random-cluster}

\def\la{\langle}
\def\ra{\rangle}
\def\eps{\epsilon}
\def\g{\gamma}

\def\Si{\Sigma}
\def\rc{random-cluster}
\def\ZZ{{\mathbb Z}}
\def\RR{{\mathbb R}}

\def\PP{{\mathbb P}}

\def\CC{{\mathbb C}}
\def\ZR{\ZZ\times\RR}
\def\ZdR{\ZZ^d\times\RR}
\def\GR{G \times\RR}

\def\VR{V\times\RR}
\def\Plb{\PP_{\lam,\de}}
\def\PLlb{\PP_{\La,\lam,\de}}

\def\PLlbq{\PP_{\La,\lam,\de,q}}
\def\PTlbq{\PP_{\lam,\de,q,T}}

\def\Plbq{\PP_{\lam,\de,q}}
\def\PLb{\PP_{\La,\de}}
\def\Om{\Omega}

\def\om{\omega}
\def\De{\Delta}

\def\th{\theta}
\def\eps{\epsilon}
\def\La{\Lambda}

\def\oo{\infty}
\def\lest{\le_{\mathrm {st}}}

\def\thetac{\theta_{\mathrm{c}}}
\def\rhoc{\rho_{\mathrm{c}}}
\def\lora{\rightarrow}

\def\be{\begin{equation}}
\def\ee{\end{equation}}
\def\sm{\setminus}
\def\resp{respectively}

\def\lra{\leftrightarrow}

\def\phmb{\phi_{m,\b}}
\def\phGb{\phi_{G,\b}}
\def\phmbp{\phmb^{\mathrm p}}
\def\phnbp{\phnb^{\mathrm p}}

\def\phnb{\phi_{n,\b}}
\def\ol#1{\overline{#1}}

\def\rad{\mathrm{rad}}

\begin{document}
\title{Space--time percolation}
\author{Geoffrey Grimmett\footnote{Centre for Mathematical Sciences, 
University of Cambridge, Wilberforce Road, Cambridge CB3 0WB, UK}} 
   \maketitle

\begin{abstract}
The contact model for the spread of disease may be viewed
as a directed percolation model on $\ZR$ in which the continuum
axis is oriented in the direction of increasing time. Techniques from
percolation have enabled a fairly complete analysis of
the contact model at and near its critical point. 
The corresponding process when the time-axis is unoriented
is an undirected percolation
model to which now standard techniques may be applied. One may
construct in similar vein a \rc\ model on $\ZR$, with associated continuum
Ising and Potts models. These models are of independent interest, 
in addition to providing a path-integral representation of the quantum Ising
model with transverse field. This representation
may be used to obtain a bound on the entanglement
of a finite set of spins in the quantum Ising model
on $\ZZ$, where this entanglement is measured via the entropy of the reduced density matrix.
The mean-field version of
the quantum Ising model gives rise to a \rc\ model on $K_n \times \RR$, thereby
extending the Erd\H os--R\'enyi random graph on the complete graph $K_n$.
\end{abstract}

\section{Introduction}\label{sec:int}
Brazil is justly famous for its beach life and
its probability community.  In harnessing the first to support
the second, a summer school of intellectual distinction
and international visibility  in probability theory
has been created. 
The high scientific stature of the 
organizers and of the wider Brazilian community
has ensured the attendance of a host of wonderful lecturers
during ten years of the Brazilian School of Probability,
and the School has attracted an international audience including many
young Brazilians who continue to leave their marks within this
crossroads subject of mathematics.
The warmth and vitality of Brazilian culture have been attractive features
of these summer schools, and invitations to participate are greatly valued.
This short review concerns two topics of recurring interest
at the School, namely percolation and the Ising model (in both its
classical and quantum forms), subject to the difference  
that one axis
of the underlying space is allowed to vary continuously.

The percolation process is arguably the most fundamental of models
for a disordered medium. Its theory is now well established, and  
several mathematics books have been written on and near
the topic, see \cite{BollRior,Cerf06,G99,Slade06}. Percolation is 
at the source of one of the most exciting areas of contemporary
probability theory, namely the theory of Schramm--L\"owner evolutions (SLE).
This theory threatens to explain the relationship between probabilistic
models and conformal field theory, and is expected to lead  ultimately
to rigorous explanations of scaling theory for a host
of two dimensional models including percolation, self-avoiding walks, 
and the Ising/Potts and \rc\ models.
See \cite{Lawler,Sch06,Smir,W02}.

Percolation theory has contributed via the \rc\ model to
the study of Ising/Potts models on a given
graph $G$, see \cite{G-RC}. The methods developed
for percolation have led also to solutions of several of the
basic questions about the contact model on $\GR$, see \cite{AJ,BG1,BG,Lig2}.
It was shown in \cite{AKN} that the quantum Ising model 
with transverse
field on $G$ may be reformulated in terms of a \rc\ model on $\GR$,
and it has been shown recently in \cite{GOS} that \rc\ arguments
may be used to study entanglement in the quantum Ising model.

In this short account of percolative processes on $\GR$ for a lattice $G$,
we shall recall in Sections \ref{sec:cp}--\ref{sec:cm}
the problems of percolation on $\GR$, and of the contact model on $G$.
This is followed in Section \ref{sec:rc} by a description of
the continuum \rc\ model on $\GR$, and its application to continuum Ising/Potts
models.  In Section \ref{sec:qim} we present 
a summary of the use of
\rc\ techniques to study entanglement in the quantum Ising model on $\ZZ$.
An account is included of a recent result of \cite{GOS} stating that
the entanglement entropy of a line of $L$ spins has order not exceeding $\log L$
in the strong-field regime. The proof relies on a property of 
\rc\ measures termed `ratio weak-mixing',  studied earlier
in \cite{Al1,Al2} for the \rc\ model on a lattice.
The corresponding mean-field model is considered in Section 6.

\section{Continuum percolation}\label{sec:cp}
Let $G=(V,E)$ be a finite or countably infinite graph which, for simplicity,
we take to be connected with neither loops nor multiple edges.
We shall usually take $G$ to be a subgraph of
the hypercubic lattice $\ZZ^d$ for some $d\ge 1$.
The models of this paper inhabit the space
$\GR$, which we refer to as space--time, and we
think of $\GR$ as being obtained by attaching a `time-line' $(-\oo,\oo)$
to each vertex $x \in V$.

Let $\lam,\de \in (0,\oo)$. The continuum percolation
model on $\GR$ is constructed via processes of
`cuts' and `bridges' as follows. For each $x\in V$, we select a Poisson
process $D_x$ of points in $\{x\}\times\RR$ with intensity $\de$; the
processes $\{D_x: x\in V\}$ are independent, 
and the points in the $D_x$ are termed
`cuts'. For each $e=\la x,y\ra \in E$, 
we select a Poisson process $B_e$ of points
in $\{e\}\times\RR$ with intensity $\lam$;
the processes $\{B_e: e\in E\}$ are independent of each other and of the
$D_x$. Let $\Plb$ denote the probability measure
associated with the family of such Poisson processes indexed by $V\cup E$.

For each $e =\la x,y\ra \in E$ and $(e,t)\in B_e$, we think of $(e,t)$
as an edge joining the endpoints $(x,t)$
and $(y,t)$, and we refer to this edge as a `bridge'.  
For $(x,s), (y,t) \in V\times\RR$, we write
$(x,s)\lra (y,t)$ if there exists a path $\pi$ with
endpoints $(x,s)$, $(y,t)$ such that: $\pi$ comprises cut-free sub-intervals
of $\GR$ together with bridges.
For $\La,\De \subseteq \VR$, we write $\La\lra \De$ if there exist
$a\in \La$ and $b\in \De$ such that $a\lra b$.

For $(x,s)\in \VR$, let $C_{x,s}$ be the set of all points
$(y,t)$ such that $(x,s)\lra(y,t)$. The clusters $C_{x,s}$ have been studied in
\cite{BG}, where the case $G=\ZZ^d$ was considered in some detail.
Let $0$ denote the origin $(0,0)\in \ZdR$, and let
$C=C_0$ denote the cluster at the origin.
Noting that $C$ is a union of line-segments, we write $|C|$ for the
Lebesgue measure of $C$. The \emph{radius} $\rad(C)$ of $C$ is given by
$$
\rad(C) = \sup\bigl\{\|x\|+|t|: (x,t)\in C\bigr\},
$$
where 
$$
\|x\|= \sup_i |x_i|, \qquad x=(x_1,x_2,\dots,x_d)\in\ZZ^d,
$$
is the supremum norm on $\ZZ^d$.

The critical point of the process is defined by
$$
\lamc(\de) = \sup\{\lam: \theta(\lam,\de) = 0\},
$$
where 
$$
\theta(\lam,\de) = \Plb(|C|=\oo).
$$
It is immediate by time-scaling that $\theta(\lam,\de)=\theta(\lam/\de,1)$,
and we shall use the abbreviations $\lamc=\lamc(1)$ and $\theta(\lam)=\theta(\lam,1)$.

The following exponential-decay theorem will be useful for the study
of the quantum Ising model in Section \ref{sec:qim}.

\begin{thm}\label{thm:cp}{\bf \cite{BG}}
Let $G=\ZZ^d$ where $d\ge 1$, and consider continuum percolation on $\GR$.
\begin{romlist}
\item 
We have that $\theta(\lamc)=0$.
\item Let $\lam,\de\in (0,\oo)$.
There exist  $\b$, $\g$ satisfying $\b,\g >0$
for $\lam/\de < \lamc$ such that{\rm:}
\begin{alignat}{2}
\Plb(|C|\ge k) &\le e^{-\g k},\qquad &&k > 0,\\
\Plb(\rad(C) \ge k) &\le e^{-\b k}, &&k > 0.
\end{alignat}
\item When $d=1$, $\lamc=1$.
\end{romlist}
\end{thm}

The situation is rather different when the environment is chosen at random.
With $G=(V,E)$ as above, suppose that the Poisson process of
cuts at a vertex $x\in V$ has some intensity $\de_x$, and that of bridges parallel
to the edge  $e= \langle x,y\rangle \in E$ has some intensity $\lam_e$.
Suppose further that the $\de_x$, $x\in V$, are independent, identically distributed
random variables, and the $\lam_e$, $e\in E$ also. Write $\De$ and $\Lambda$
for independent random variables having the respective distributions, and $P$
for the probability
measure governing the environment.
[As before, $\Plb$ denotes the measure associated with the percolation model in
the given environment.]

If there exist $\lam',\de'\in (0,\oo)$ such that $\lam'/\de' < \lamc$
and $P(\Lambda \le \lam')=P(\Delta \ge \de') = 1$, then
the process is almost surely dominated by a subcritical 
percolation process, whence there is (almost sure) exponential decay in the sense of
Theorem \ref{thm:cp}(ii).  This may fail in an interesting
way if there is no such almost-sure domination, in that one may prove exponential decay
in the space-direction but only a weaker decay in the time-direction.

For any probability measure $\mu$ and function $f$,
we write $\mu(f)$ for the expectation of $f$ under $\mu$.
For $(x,s), (y,t) \in \ZZ^d\times \RR$ and $q\ge 1$, we define
$$
d_q(x,s;y,t) = \max\bigl\{\|x-y\|, [\log(1+|s-t|)]^q\bigr\}.
$$

\begin{thm}\label{thm:klein}{\bf\cite{Klein1,Klein2}}
Let $G=\ZZ^d$ where $d \ge 1$. Suppose
that
$$
\Gamma = \max\left\{P\bigl([\log(1+\Lambda)]^\b\bigr), 
P\bigl( [\log(1+\Delta^{-1})]^\b\bigr) \right\} <\oo,
$$
for some $\b>2d^2\bigl(1+\sqrt{1+d^{-1}} + (2d)^{-1}\bigr)$.
There exists $Q=Q(d,\b)>1$ such that the following holds. For $q\in [1,Q)$ and $m>0$,
there exists $\eps=\eps(d,\b,\Gamma,m,q)>0$ and $\eta=\eta(d,\b,q)>0$ such that{\rm:} if
$$
E\Bigl(\bigl[\log(1+(\Lambda/\Delta))\bigr]^\b\Bigr) < \eps,
$$
there exist identically distributed random variables 
$D_x\in L^\eta(P)$, $x\in\ZZ^d$, such that
$$
\Plb\bigl((x,s) \lra (y,t)\bigr)\le 
\exp\bigl[-m d_q(x,s;y,t)\bigr] \qquad \text{if } d_q(x,s;y,t)\ge D_x,
$$
for $(x,s),(y,t) \in \ZZ^d\times\RR$. 
\end{thm}  

The corresponding theorem of \cite{Klein1}
contains no estimate for the tail of the $D_x$. The above moment property
may be derived from the Borel--Cantelli argument used in the proof
of \cite{Klein1}, which proceeds by a so-called multiscale analysis,
see \cite{GOS}, Section 8. 
Explicit values may be given for the constants $Q$ and $\eta$, namely
$$
Q = \frac{\b(\a-d+\a d)}{\a d(\a+\b+1)},
$$
where $\a = d+\sqrt{d^2+d}$, and
one may take $\eta$ satisfying
$$
(\eta+1)\a < \frac \b\a\left(\frac{\a-d+\a d}q -\a d\right) -d.
$$

Complementary
accounts of the \emph{survival} of the process in a random environment
may be found in \cite{AKN,Andjel,CKP,NewmanV}.

We mention two further types of `continuum' percolation processes
that arise in applications and have attracted the attention of probabilists. 
Let
$\Pi$ be a Poisson process of points in $\RR^d$ with
intensity $1$. Two points $x,y\in\Pi$ are joined by an edge, 
and said to be adjacent,
if they satisfy a given condition of proximity. One now asks
for conditions under which the resulting random graph possesses
an unbounded component.

The following conditions of proximity have been studied in the literature.

\begin{numlist}
\item \emph{Lily-pond model.} Fix $r > 0$, and join $x$ and $y$ if and only if 
$|x-y|\le r$, where $|\cdot|$ denotes Euclidean distance. There has
been extensive study of this process, 
and of its generalization, the \emph{random
connection model}, in which $x$ and $y$ are joined with probability
$g(|x-y|)$ for some given non-increasing function $g:(0,\oo)\to[0,1]$.
See \cite{G99,MeeR,Pen03}.
 
\item \emph{Voronoi percolation.} To each $x\in\Pi$ we associate the \emph{tile}
$$
T_x=\bigl\{z\in \RR^d: |z-x| \le |z-y|\mbox{ for all } y\in\Pi\sm\{x\}\bigr\}.
$$
Two tiles $T_x$, $T_y$ are declared \emph{adjacent} if their boundaries
share a facet of a hyperplane of $\RR^d$.
We colour each tile \emph{red} with probability $\rho$, 
different tiles receiving
independent colours, and we ask for conditions under which there
exists an infinite path of red tiles.

This model has a certain property of conformal invariance when $d=2,3$, see \cite{BenS98}.
When $d=2$, there is an obvious property of self-matching, leading to
the conjecture that the critical point is given by
$\rhoc=\frac12$, and this has been
proved recently in \cite{BollRior2}.
\end{numlist}

\section{The contact model}\label{sec:cm}
Just as directed percolation on $\ZZ^d$ arises by allowing only open paths that are `stiff'
in one direction, so the contact model on $G$ is obtained from percolation on $\GR$
by requiring that open paths traverse time-lines in the direction of \emph{increasing} time.

As before, we let $D_x$, $x\in V$, be Poisson processes with intensity
$\de$, and we term points in the $D_x$ `cuts'. We replace each $e=\la x,y\ra\in E$ by
two oriented edges $[x,y\ra$, $[y,x\ra$, the first oriented from  $x$ to $y$, and the 
second from $y$ to $x$. Write $\vec E$ for the set of oriented edges thus obtained from $E$.
For each $\vec e = [x,y\ra \in \vec E$, we let $B_{\vec e}$ be a Poisson process with
intensity $\lam$; members of $B_{\vec e}$ are termed `directed bridges'
from $x$ to $y$.

For $(x,s),(y,t) \in \VR$, we write $(x,s)\lora (y,t)$ if
there exists an oriented  path $\pi$ from
$(x,s)$ to $(y,t)$ such that: $\pi$ comprises cut-free sub-intervals
of $\VR$ traversed in the direction of increasing time, 
together with directed bridges
in the directions of their orientations.
For $\La,\De \subseteq \VR$, we write $\La\lora \De$ if there exist
$a\in \La$ and $b\in \De$ such that $a\lora b$.

The directed cluster $\vec C$ at the origin is the set 
$$
\vec C = \bigl\{(x,s)\in \VR: 0 \lora (x,s)\bigr\},
$$
of points reachable from the origin $0$ along paths directed away from $0$.
The percolation probability is given by
$$
\vec\theta(\lambda,\de) = \Plb(|\vec C|=\oo),
$$
and the critical point by
$$
\vlamc(\de) = \sup\bigl\{\lam:\vec\theta(\lam,\de)=0\bigr\}.
$$
As before, we write $\vec\theta(\lam) = \vec\theta(\lam,1)$ and $\vlamc=\vlamc(1)$.

Parts (i) and (ii) of Theorem \ref{thm:cp} are valid in this new setting,
with $C$ replaced by $\vec C$, etc, see \cite{BG1}. The exact value
of the critical point is unknown even when $d=1$, 
although there are physical reasons to believe in this case that
$\vec\lamc=1.694\dots$, the critical value of the so-called reggeon
spin model, see \cite{GT,Lig1}. The contact model in a random environment
may be studied as in Theorem \ref{thm:klein}.

Further theory of the contact model may be found in \cite{Lig1,Lig2}.
Sakai and Slade \cite{SS} have shown how to apply the lace expansion to the
spread-out contact model on $\ZZ^d$ for $d > 4$, and related results
are valid for oriented percolation even when the connection function
has unbounded domain, see \cite{CSa,CSh}.

\section{Random-cluster and Ising/Potts models}\label{sec:rc}
The percolation model on a graph $G=(V,E)$ 
may be generalized to obtain the
\rc\ model on $G$, see \cite{G-RC}. 
Similarly, the continuum percolation model
on $\GR$ may be extended to a continuum \rc\ model. 
Let $W$ be a finite subset of $V$ that
induces a connected subgraph of $G$, and let $E_W$ denote the
set of edges joining vertices in $W$. Let $T\in(0,\oo)$, and let
$\La$ be the `box' $\La=W \times [0,T]$. Let $\PLlb$ denote
the probability measure associated with the Poisson processes
$D_x$, $x\in W$, and $B_e$, $e=\la x,y\ra \in E_W$.
As sample space we take the
set $\Om_\La$ comprising all finite sets of cuts and bridges in $\La$,
and we may assume without loss of generality that no cut is the endpoint of any bridge.
For $\om\in\Om_\La$, we write $B(\om)$ and $D(\om)$ for the sets
of bridges and cuts, respectively, of $\om$. The appropriate $\s$-field 
$\sF_\La$ is that generated by the open 
sets in the associated Skorohod topology,
see \cite{BG,EK}. 
 
For a given configuration $\om\in\Om_\La$, let $k(\om)$
be the number of its clusters under the connection relation $\lra$. 
Let $q\in(0,\oo)$, and define the
`continuum \rc' probability measure
$\PLlbq$ by
\be\label{rcPo}
d\PLlbq(\om) = \frac 1Z q^{k(\om)}d\PLlb(\om),
\qq \om\in\Om_\La,
\ee
for an appropriate normalizing constant, or `partition function', $Z=Z_\La(\lam,\de,q)$.
The quantity $q$ is called the \emph{cluster-weighting factor}.
The continuum \rc\ model may be studied in very much
the same way as the \rc\ model on a lattice, see \cite{G-RC}.

The space $\Om_\La$ is a partially ordered space with 
order relation given by: $\om_1\le\om_2$ if $B(\om_1)\subseteq B(\om_2)$ and
$D(\om_1) \supseteq D(\om_2)$. A random variable $X:\Om_\La\to\RR$
is called \emph{increasing}
if $X(\om)\le X(\om')$ whenever $\om\le\om'$. An event $A\in\sF_\La$ is
called \emph{increasing} if its indicator function $1_A$ is increasing.
Given two probability measures $\mu_1$, $\mu_2$ on
a measurable pair $(\Om_\La,\sF_\La)$, we write
$\mu_1\lest \mu_2$ if $\mu_1(X) \le \mu_2(X)$
for all bounded increasing continuous random variables $X:\Om_\La\to\RR$.

The measures $\PLlbq$ have certain
properties of stochastic ordering as the parameters $\La$, $\lam$,
$\de$, $q$ vary. The basic theory will be assumed here, and the reader is
referred to  \cite{Bj}
for further details. In rough terms, the $\PLlbq$ inherit the properties
of stochastic ordering and positive association enjoyed by their counterparts 
on discrete graphs. Of particular value later will
be the stochastic inequality
\begin{equation}\label{stochcomp}
\PLlbq \lest \PLlb\qquad\mbox{when } q \ge 1.
\end{equation}

While it will not be important for what follows, we note that
the thermodynamic limit may be taken in much the same manner
as for the discrete \rc\ model, whenever $q \ge 1$. Suppose, for example, that $W$ is
a finite connected subgraph of the lattice  $G=\ZZ^d$, and assign
to the box $\La = W \times [0,T]$ a suitable boundary condition. 
As in \cite{G-RC},
if the boundary condition $\tau$ is chosen in such a way that the measures
$\PLlbq^\tau$ are monotonic as $W\uparrow \ZZ^d$, then the weak limit
$\PTlbq^\tau = \lim_{W\uparrow\ZZ^d} \PLlbq^\tau$ exists. One may
similarly allow the limit as $T\to\oo$ to obtain a measure
$\Plbq^\tau=\lim_{T\to\oo} \PTlbq^\tau$. 

Let $G=\ZZ^d$.
Restricting ourselves for convenience to the case
of free boundary conditions, we define the percolation probability by
$$
\theta(\lam,\de,q) = \Plbq(|C_0| = \oo),
$$
and the critical point by
$$
\lamc(\ZZ^d, q)= \sup\bigl\{\lam: \theta(\lam,1,q) = 0\bigr\}.
$$
In the special case $d=1$, the \rc\ model has a property of self-duality
that leads to the following conjecture.

\begin{conj}\label{conj:critpt}
The continuum \rc\ model on $\ZR$ with $q \ge 1$ has critical value
$\lamc(\ZZ,q) = q$.
\end{conj}
 
It may be proved by standard means that $\lamc(\ZZ,q) \ge q$. 
See \cite{G-RC}, Section 6.2, for the corresponding result on the discrete lattice $\ZZ^2$.

The \emph{continuum Potts model} on $\GR$ is given as follows.
Let $q\in\{2,3,\dots\}$. To each cluster of the \rc\ model with
cluster-weighting factor $q$ is assigned a `spin' from the space
$\Si = \{1,2,\dots,q\}$, different clusters receiving independent spins.
The outcome is a function $\s:\VR \to \Si$, and this is the spin-vector
of a `continuum $q$-state Potts model' with parameters $\lam$ and $\de$.
When $q=2$, we refer to the model as a continuum Ising model.

It may be seen that the law of the above spin model on
$\La = W\times[0,T]$ is given by
$$
d\PP(\s) = \frac 1Z e^{\lam L(\s)}\,d\PLb(D_\s) ,
$$
where $D_\s$ is the set of $(x,s)\in W\times[0,T]$ such that
$\s(x,s-) \ne \s(x,s+)$, $\PLb$ is the law of a family
of independent Poisson processes on the time-lines $\{x\}\times[0,T]$,
$x\in W$, with intensity $\de$, and 
$$
L(\s) =  \sum_{\la x,y\ra \in E_W} \int_0^T 1_{\{\s(x,u)=\s(y,u)\}}\,du
$$
is the aggregate Lebesgue measure of those subsets of adjacent time-lines
on which the spins are equal. As usual, $Z$ is an appropriate constant.

The continuum Ising model has arisen in the study by 
Aizenman, Klein, and Newman, \cite{AKN}, of the quantum
Ising model with transverse field, as described in the next section.

\section{The quantum Ising model}\label{sec:qim}

Aizenman, Klein, and Newman reported in \cite{AKN} a representation of
the quantum Ising model in terms of the $q=2$ continuum \rc\ and Ising models. 
This was motivated in part by \cite{CKP} and by earlier work
referred to therein. We summarise
this here, and we indicate how it may be used to study the property of entanglement 
in the quantum Ising model on $\ZZ$.

The quantum Ising model on a finite graph $G=(V,E)$ is given as follows.
To each vertex $x\in V$ is associated a quantum spin-$\frac12$ with local Hilbert space
$\CC^2$. The Hilbert space $\mathcal{H}$ for the system is therefore
the tensor product $\mathcal{H} = \bigotimes_{x\in V} \CC^2$. As basis
for the copy of $\CC^2$ labelled by $x$, we take the two eigenstates, 
denoted as
$|+\rangle_x = \left(\begin{matrix}1\\0\end{matrix}\right) $  and
$|-\rangle _x= \left(\begin{matrix}0\\1\end{matrix}\right) $, 
of the Pauli operator
$$
\sigma^{(3)}_x = \left(
\begin{array}{cc} 1 & 0
\\ 0 & -1\end{array} \right)
$$
at the site $x$, with corresponding eigenvalues $\pm 1$.
The other two Pauli operators with respect
to this basis are the matrices
\begin{equation}
\sigma^{(1)}_x= \left( \begin{array}{cc} 0 & 1 \\ 1 & 0\end{array}
\right), \qquad \sigma^{(2)}_x= \left ( \begin{array}{cc} 0& -i \\
i & 0\end{array}\right).
\end{equation}
In the following, $|\phi\rangle$ denotes a vector and $\langle \phi|$
its adjoint. 

Let $D$ be the set of $2^{|V|}$
basis vectors $|\eta\rangle$ for $\sH$ of the form $|\eta\rangle =
\bigotimes_x|\pm\rangle_x$. There is a natural one--one correspondence between $D$
and the space $\Si= \Si_V = \prod_{x\in V}\{-1,+1\}$.
We shall sometimes speak of members of $\Si$ as basis vectors,
and of $\sH$ as the Hilbert space
generated by $\Si$. 

The Hamiltonian of the quantum Ising model
with transverse field is the operator
\begin{equation}\label{ham}
H = -\tfrac{1}{2} \lam \sum_{e=\langle x,
y\rangle\in E} \sigma^{(3)}_x
 \sigma^{(3)}_y - \delta \sum_{x\in V} \sigma^{(1)}_x,
\end{equation}
generating the operator $e^{-\b H}$ where $\b$ denotes inverse
temperature. Here,  $\lambda, \de\geq 0$  are
the spin-coupling and external-field intensities,
respectively. The 
Hamiltonian has a unique pure ground state $|\psi_G \rangle$
defined at zero-temperature (that is, in the limit as $\b\to\oo$) 
as the eigenvector corresponding to
the lowest eigenvalue of $H$.

Let
\be\label{eq:operb}
\rho_G(\b) = \frac 1 {Z_G(\b)}e^{-\b H},
\ee
where
$$
Z_G(\b) = \tr (e^{-\b H})
= \sum_{\eta\in \Si} \la\eta| e^{-\b H} |\eta\ra.
$$
It turns out that the matrix elements of $\rho_G(\b)$ may be
expressed as a type of `path integral' with respect to the continuum \rc\ model
on $G\times [0,\b]$ with parameters $\lam$, $\de$ and $q=2$. Let 
$\La=V\times[0,\b]$, write $\Om_\La$
for the configuration space of the latter model, and let $\phGb$ be the appropriate
continuum \rc\ measure on $\Om_\La$ (with free boundary conditions). 
For $\om\in\Om_\La$, let $S_\om$ denote
the space of all functions $s: V\times[0,\b] \to \{-1,+1\}$ that
are constant on the clusters of $\om$, and let $S$ be the union of the
$S_\om$ over $\om\in\Om_\La$. Given  $\om$, we may pick an element
of $S_\om$ uniformly at random, and we denote this random
element as $\s$. We shall
abuse notation by using $\phGb$ to denote the 
ensuing probability measure on the
coupled space $\Om_\La\times S$. For $s\in S$ and $W\subseteq V$, we write
$s_{W,0}$ (\resp, $s_{W,\b}$) for the vector $(s(x,0): x\in W)$
(\resp, $(s(x,\b): x\in W)$). We abbreviate $s_{V,0}$ and $s_{V,\b}$ to
$s_0$ and $s_\b$, \resp.

The following representation of the matrix elements of $\rho_G(\b)$
is obtained by expanding the exponential in \eqref{eq:operb}, and it
permits the use of \rc\ methods to study the matrix $\rho_G(\b)$.
For example, as pointed out in \cite{AKN}, it implies the
existence of the low-temperature limits
$$
\langle \eta'|\rho_G|\eta \rangle =
\lim_{\b\to\oo} \langle \eta'|\rho_G(\b)|\eta \rangle,
\qquad \eta,\eta'\in\Si.
$$

\begin{thm}\label{thm:akn} {\bf\cite{AKN}}
The elements of the
density matrix $\rho_G(\b)$ are given by
\be
\langle \eta'|\rho_G(\b)|\eta \rangle = \frac
{\phGb(\s_0=\eta,\, \s_\b=\eta')}{\phGb(\s_0=\s_\b)},
\qquad \eta,\eta'\in \Si.
\label{eq:thing}
\ee
\end{thm}

This representation may be used to study entanglement in the quantum
Ising model on $G$. Let $W \subseteq V$, and consider the
\emph{reduced density matrix} 
\be\label{reddo} 
\rho_G^W(\b) =
 \tr_{V\setminus W} (\rho_G(\b)),
\ee
where the trace is performed over the
Hilbert space $\sH_{V\sm W} = \bigotimes_{x\in V\sm W}
\CC^2$ of the spins belonging to $V\sm W$. By an analysis
parallel to that leading to Theorem \ref{thm:akn}, we obtain
the following.

\begin{thm}\label{thm:gos} {\bf\cite{GOS}}
The elements of the reduced
density matrix $\rho_G^W(\b)$ are given by
\be
\langle \eta'|\rho_G^W(\b)|\eta \rangle = \frac
{\phGb(\s_{W,0}=\eta,\, \s_{W,\b}=\eta'\mid E)}{\phGb(\s_0=\s_\b\mid E)},
\qquad \eta,\eta'\in \Si_W,
\label{eq:thing2}
\ee
where $E$ is the event that $\s_{V\sm W,0}=\s_{V\sm W,\b}$.
\end{thm}

Let $D_W$ be the set of $2^{|W|}$
vectors $|\eta\rangle$ of the form $|\eta\rangle =
\bigotimes_{x\in W}|\pm\rangle_x$, and write $\sH_W$ for
the space generated by $D_W$. Just as before, 
there is a natural one--one correspondence between $D_W$
and the space $\Si_W=\prod_{x\in W}\{-1,+1\}$,
and we shall regard $\sH_W$ as the Hilbert space
generated by $\Si_W$.

We may write
$$
\rho_G=\lim_{\b\to\oo}\rho_G(\b) =|\psi_G \rangle\langle\psi_G|
$$
for the density matrix corresponding to the ground state of the
system, and similarly
\be\label{reddo2} 
\rho_G^W =
\tr _{V\setminus W} (|\psi_G\rangle\langle \psi_G|)
= \lim_{\b\to\oo}\rho_G^W(\b) . 
\ee

There has been extensive study of entanglement in the physics literature,
see the references in \cite{GOS}. The entanglement of the spins
in $W$ may be defined as 
follows.

\begin{defn}\label{ent}
The \emph{entanglement} of the vertex-set $W$ 
relative to its complement $V \setminus W$ is the entropy
\begin{equation}\label{entdef}
S_G^W= -\tr(\rho_G^W \log_2 \rho_G^W).
\end{equation}
\end{defn}

The behaviour of $S_G^W$, for general $G$ and $W$, is not
understood at present. Instead, we specialise here to the case
of a finite subset of the one-dimensional lattice $\ZZ$.
Let $m,L \ge 0$ and take $V=[-m, m+L]$ and $W=[0,L]$,
viewed as subsets of $\ZZ$.
We obtain $G$ from $V$ by adding edges between each pair
$x,y\in V$ with $|x-y|=1$.
We write $\rho_m(\b)$ for $\rho_G(\b)$, and $S_m^L$ for
$S_G^W$.
A key step in the study of $S_m^L$ for large $m$ is a bound
on the norm of the difference $\rho_m^L-\rho_n^L$. For an operator
$A$ on $\sH$, let
$$
\|A\| = \sum_{\|\psi\|=1} \bigl| \la \psi| A |\psi\ra \bigr|,
$$
where the supremum is over all $\psi\in \sH_L$
with $L^2$-norm $1$.

\begin{thm} {\bf\cite{GOS}}\label{thm:expdecrc}
Let $\lam,\de\in (0,\oo)$ and write $\th=\lam/\de$. 
There exist constants $C$, $\a$, $\g$ depending
on $\th$ and satisfying $\g > 0$ when $\th < 1$  such that{\rm:}
\be\label{eq:normdiff}
\| \rho_m^L-\rho_n^L\| \le \min\bigl\{2, CL^\a e^{-\g m}\bigr\},\qquad
2\le m\le n<\oo.
\ee
\end{thm}

One would
expect that $\g$ may be taken in such a manner that
$\g>0$ under the weaker assumption $\lam/\de<2$, but
this has not yet been proved (cf.\ Conjecture \ref{conj:critpt}). 

Inequality \eqref{eq:normdiff}
is proved in \cite{GOS} by the following route. 
Consider the \rc\ model with $q=2$ on the space--time
graph $\La = V\times [0,\b]$ with `partial periodic top/bottom boundary
conditions'; that is, for each $x\in V\sm W$, 
we identify the two vertices
$(x,0)$ and $(x,\b)$. Let $\phmbp$ denote the
associated \rc\ measure on $\Om_\La$. To each cluster
of $\om$ ($\in\Om_\La$) we assign a random spin
from $\{-1,+1\}$ in the usual manner, and we
abuse notation by using $\phmbp$ for the measure
governing both the \rc\ configuration and the spin configuration.
Let $\amb=\phmbp(\s_{W,0}=\s_{W,\b})$, noting that
$\amb = \phmb(\s_{0}=\s_{\b}\mid E)$ as in \eqref{eq:thing2}.

By Theorem \ref{thm:gos},
\be\label{eq:4}
\langle \psi| \rho_m^L(\b)-\rho_n^L(\b) |\psi\rangle
= \frac{\phmbp(c(\s_{W,0})\ol{c(\s_{W,\b})})}{\amb}
- \frac{\phnbp(c(\s_{W,0})\ol{c(\s_{W,\b)}})}{\anb},
\ee
where $c:\{-1,+1\}^W \to \CC$ and
$$
\psi=\sum_{\eta\in \Si_W} c(\eta) \eta \in \sH_W.
$$

The \rc\ property of ratio weak-mixing
is used in the derivation of \eqref{eq:normdiff} from
\eqref{eq:4}. At the final step of the proof of Theorem \ref{thm:expdecrc}, the \rc\ model is compared
with the continuum percolation
model of Section \ref{sec:cp}, and the exponential decay 
of Theorem \ref{thm:expdecrc} follows
by Theorem \ref{thm:cp}. A logarithmic bound on the entanglement entropy follows
for sufficiently small $\lam/\de$.

\begin{thm}\label{entest} {\bf\cite{GOS}}
Let $\lam,\de\in(0,\oo)$ and write $\th=\lam/\de$. 
There exists $\theta_0\in(0,\oo)$ such that{\rm:}
for $\th <\theta_0$, there exists $K = K(\th) <\oo$ such that
\begin{equation}
S_m^L \le K \log_2 L,\qquad m\ge 0,\ L \ge 2.
\end{equation}
\end{thm}

A stronger result is expected, namely that the 
entanglement $S_m^L$ is bounded above, uniformly in $L$,
whenever $\th$ is sufficiently small, and perhaps
for all $\th<\thetac$ where $\thetac=2$ is the critical point.
See Conjecture \ref{conj:critpt} and the references in \cite{GOS}.
There is no rigorous picture known of the behaviour
of $S_m^L$ for large $\th$, or of the corresponding
quantity in dimensions $d\ge 2$, although Theorem \ref{thm:expdecrc} has
a counterpart in this setting. Theorem \ref{entest} may be extended
to the disordered system in which the intensities $\lam$, $\de$ 
are independent random variables indexed by the vertices and
edges of the underlying graph, subject to certain
conditions on these variables (cf.\  Theorem \ref{thm:klein}
and the preceding discussion).

\section{The mean-field continuum model}
The term `mean-field' is often interpreted in percolation theory as
percolation on either a tree (see \cite{G99}, Chapter 10) or a complete graph.
The latter case is known as the
Erd\H os--R\'enyi random graph $G_{n,p}$, and this is the random
graph obtained from the complete graph $K_n$ on $n$ vertices by deleting
each edge with
probability $1-p$. The theory of $G_{n,p}$ is well developed
and rather refined, see \cite{Boll01,JLR}, and particular
attention has been paid to the emergence of the giant cluster
for $p=\lam/n$ and $\lam\simeq 1$. A similar theory
has been developed for the \rc\ model on $K_n$
with parameters $p$, $q$, see \cite{BGJ,G-RC,LL}.

Unless boundary conditions are introduced
in the manner of \cite{GJ0,Hagg}, the continuum \rc\ model on
a tree may be solved exactly by standard means.
We therefore concentrate here on the case of the complete graph
$K_n$ on $n$ vertices.
Let $\b>0$, and attach to each vertex the line $[0,\b]$ with its endpoints
identified; thus, the line forms a circle. We now consider the 
continuum \rc\ model
on $K_n\times[0,\b]$ with parameters $p=\lam/n$, $\de=1$, and $q$. 
[The convention of setting $\de=1$ differs from that of \cite{IoffeL} 
but is consistent with  that adopted in earlier work on related models.]

Suppose 
that $q \ge 1$, so that we may use methods based on stochastic
comparisons. It
is natural to ask for the critical value $\lamc =\lamc(\b,q)$ of $\lam$ above
which the model possesses a giant cluster. This has been
answered by Ioffe and Levit, \cite{IoffeL},
in the special case $q=1$.
Let $F(\b,\lam)$ be given by
$$
F(\b,\lam) = \lam\bigl[2(1-e^{-\b})-\b e^{-\b}\bigr],
$$
and let $\lamc = \lamc(\b)$ be chosen so that $F(\b,\lamc)=1$.

\begin{thm} {\bf\cite{IoffeL}}\label{IL}
Let $M$ be the maximal
(one-dimensional) Lebesgue measure
of the clusters of the process with parameters $\b$, $p=\lam/n$, $\de=1$, $q=1$.
Then, as $n\to\oo$,
$$
\frac 1n M \to\begin{cases} 0 &\mbox{if } \lam < \lamc,\\
\b \pi &\mbox{if } \lam > \lamc,
\end{cases}
$$
where $\pi=\pi(\b,\lam)\in(0,1)$ when $\lam> \lamc$, and the convergence
is in probability.
\end{thm}

When $\lam > \lamc$, the density of the giant cluster
is $\pi$, in that there is probability $\pi$ that any
given point of $K_n\times[0,\b]$ lies in this 
giant cluster.
The proof of Theorem \ref{IL} is simple to motivate. 
Let $0$ be a vertex of $K_n$, and let $I$
be the maximal cut-free interval of $0\times[0,\b]$
(viewed as a circle) containing the point $0\times 0$. Given $I$,
the mean number of
bridges leaving $I$ is $\lam|I|(n-1)/n \sim \lam |I|$, where $|I|$ is the Lebesgue
measure of $I$. One may thus approximate to the cluster
at $0\times 0$ by a branching process with mean family-size
$\lam E|I|$. It is elementary that $\lam E|I| = F(\b,\lam)$, 
which is to say that the branching process is subcritical
(\resp, supercritical) if $\lam < \lamc$ (\resp, $\lam > \lamc$). 
The details of the proof may be found in \cite{IoffeL}, and
a further proof has appeared in \cite{Janson06}. The quantity $\pi$
is of course the survival probability
of the above branching process, and this may be calculated
in the standard way on noting that $|I|$
is distributed as $\min\{U+V,\b\}$ where
$U$, $V$ are independent, exponentially distributed, random variables
with mean 1.

What is the analogue of Theorem \ref{IL} when $q\ne 1$? Indications
are presented in \cite{IoffeL} of the critical value when $q=2$, and the problem
is posed there of proving this value by calculations of the \rc\ type to be 
found in \cite{BGJ}. There is a simple argument that yields 
upper and lower bounds for
the critical value for any $q\in[1,\oo)$. We present this next, and also
explain our reason for believing the upper bound
to be exact when $q\in [1,2]$.

Consider the continuum \rc\ model on $K_n\times [0,\b]$ with
parameters $p=\lam/n$, $\de=1$, and $q\in(0,\oo)$. Let
\be\label{eq:ei}
F_q(\b,\lam) = \frac{\lam}{q^2}\cdot\frac{2e^{\b q}-2+\b q(q-2)}
{e^{\b q} +q-1},
\ee
noting that $F_1=F$.
 
\begin{thm} \label{thm:grg}
Let $M_q$ be the maximal
(one-dimensional) Lebesgue measure
of the clusters of the process with parameters $\b$, $p=\lam/n$, $\de=1$, 
$q \in [1,\oo)$.
\begin{romlist}
\item 
We have that
$\lim_{n\to\oo}n^{-1}M_q = 0$ if $F_q < q^{-1}$, where the convergence
is in probability.
\item
There exists $\pi_q=\pi_q(\b,\lam)$, satisfying
$\pi_q >0$ whenever $F_q>1$,
such that
$$
\liminf_{n\to\oo} P\left(\frac 1n M_q \ge \b\pi_q\right) \to 1.
$$
\end{romlist}
\end{thm}

The bound $\pi_q$ may be calculated by a branching-process argument,
in the same manner as was $\pi=\pi_1$,
above. We conjecture that $n^{-1}M_q\to 0$ in probability if $F_q<1$ and
$q\in[1,2]$. This conjecture is motivated
by the evidence of \cite{BGJ} that, in the second-order phase transition
occurring when $q\in [1,2]$, the location of the critical point
is given by the branching-process approximation
described in the sketch proof below. This amounts to the claim that the critical value
$\lamc(q)$ of the continuum \rc\ model with cluster-weighting factor $q$ satisfies
\be\label{eq:exact}
\lamc(q) = {q^2}  
\frac {e^{\b q} +q-1}{2e^{\b q}-2+\b q(q-2)},\qquad
q \in[1,2].
\ee
This is implied by Theorem \ref{IL} when $q=1$,
and by the claim of \cite{IoffeL} when $q=2$. 
Note the relatively simple formula when $q=2$,
\be\label{eq:meanf}
\lamc(2) = \frac 2{\tanh \b},
\ee
which might be termed the critical point of the \emph{quantum random graph}.
Dmitry Ioffe has pointed out that the exact calculation \eqref{eq:meanf}
may be derived from the results of \cite{Dor,FSV}. Results similar to
those of Theorem \ref{thm:grg} may be obtained for $q<1$ also.

\begin{proof}[Sketch proof of Theorem \ref{thm:grg}]
We begin with part (ii).
The idea is to bound the process below by a random graph to which
the results of \cite{BJR,Janson06} may be applied directly. The
bounding process is obtained as follows. First, we place the cuts
on each of the time-lines $x\times[0,\b]$, and we place
no bridges. Thus, the cuts on a given time-line are placed
in the manner of the continuum \rc\ model on that line.
It may be seen that the number $D$ of cuts on
any given time-line has mass function
$$
P(D=k) = \frac {e^{-\b}}Z\cdot \frac{q^{k\vee 1}\b^k}{k!},\qquad k \ge 0,
$$
where $a \vee b =\max\{a,b\}$, and $Z$ is the requisite constant,
$$
Z=(q-1)e^{-\b} + e^{\b(q-1)}.
$$
It is an easy calculation that the maximal cut-free interval $I$ containing
the point $0\times 0$ satisfies $E|I|=qF_q/\lam$.

We next place edges between pairs of time-lines according to
independent Poisson processes with intensity $\lam/ q$. We term
the ensuing graph a `product \rc\ model', and we claim that this model
is dominated (stochastically) by the continuum \rc\ model.
This may be seen in either of two ways: one may 
apply suitable comparison inequalities (see \cite{G-RC}, Section 3.4) to
a discrete approximation of $K_n\times[0,\b]$
and then pass to the continuum limit, or one may establish
it directly for the continuum model. Related material has appeared
in \cite{GeoT,Preston}. 

If this `product' \rc\ model possesses a giant cluster, then so does the original \rc\ model.
The former model may be studied either via the general techniques of
\cite{BJR,Janson06} for inhomogeneous random graphs, or
using the usual branching process approximation. We follow the latter route
here.   In the limit as $n\to\oo$,
the mean number of offspring of $0\times 0$ approaches
$(\lam/q)E|I| = F_q$, so that the branching process is supercritical
if $F_q>1$. The claim of part (ii) follows.

For part (i) one proceeds similarly, but with $\lam/q$ replaced by $\lam$
and the domination reversed. 
\end{proof}

\section{Acknowledgements} The author thanks Carol Bezuidenhout for encouraging
him to persevere with continuum percolation and the contact model many years ago,
and to Tobias Osborne and Petra Scudo for explaining the relationship between the quantum
Ising model and the continuum \rc\ model. He is grateful to Svante Janson
for their discussions of \rc\ models on complete graphs. Dima Ioffe has kindly
pointed out the link between the quantum random graph and the work
of \cite{Dor,FSV}.

\bibliography{cperc}
\bibliographystyle{plain}

\end{document}